\documentclass[11pt, english]{article}
\font\smallit=cmti10

\usepackage[english]{babel}
\usepackage{amssymb,amsthm,amsmath,amsfonts,mathrsfs}
\usepackage{pst-all,pst-3dplot,pstricks,pstricks-add,pst-math,pst-xkey}
\usepackage{graphicx}
\usepackage{latexsym}
\usepackage{ifthen, bigstrut}
\usepackage{accents}
\usepackage[english]{babel}
\usepackage{hyperref}
\usepackage{breakurl}
\usepackage{stmaryrd}
\usepackage{ifthen, bigstrut}
\usepackage{enumerate}
\usepackage{epsfig}
\usepackage{subfigure}
\usepackage{skak}
\usepackage{chessboard}
\usepackage{chessfss}
\usepackage{pst-all}
\usepackage{lscape}
\usepackage{courier}
\usepackage[all,cmtip]{xy}
\usepackage{bbold}

\newtheorem{theorem}{Theorem}

\newtheorem{corollary}[theorem]{Corollary}
\newtheorem{lemma}[theorem]{Lemma}
\usepackage{mathrsfs}
\newcommand{\Ne}{\mathscr{N}}
\renewcommand{\Pr}{\mathscr{P}}
\newcommand{\Le}{\mathscr{L}}
\newcommand{\Ri}{\mathscr{R}}

\newcommand{\A}{\mathcal{A}}

\newcommand{\Np}{\mathbb{Np}}
\newcommand{\E}{\mathbb{E}}

\newcommand{\D}{\mathbb{D}}

\newcommand{\U}{\mathbb{U}}
\newcommand{\GL}{{G^\mathcal{L}}}
\newcommand{\GR}{{G^\mathcal{R}}}
\newcommand{\HL}{H^\mathcal{L}}
\newcommand{\HR}{H^\mathcal{R}}

\usepackage{mathabx}
\newcommand{\conj}[1]{\overset{\curvearrowleftright}{#1}}

\newcommand{\pura}[2]{\{\emptyset^{#1}\!\mid \!\emptyset^{#2}\}}
\newcommand{\cg}[2]{\left\{ #1\!\mid \!#2\right\}}
\newcommand{\cgs}[2]{\{ #1\!\mid \!#2\}}
\newcommand{\atom}[1]{\emptyset^{#1}}
\renewcommand{\hat}{\widehat}

\renewcommand{\ge}{\geqslant}
\renewcommand{\le}{\leqslant}
\newcommand{\su}{\succcurlyeq}
\newcommand{\pr}{\preccurlyeq}

\newcommand{\bs}[1]{\boldsymbol#1}
\newcommand{\zs}[1]{\widetilde{#1}}
\newcommand{\abss}{absolute}

\newcommand{\abs}{\abss\ }

\renewcommand\emptyset \varnothing

\theoremstyle{definition}
\newtheorem{definition}[theorem]{Definition}

\newtheorem{problem}[theorem]{Problem}
\newtheorem{observation}[theorem]{Observation}

\emergencystretch=\maxdimen

\begin{document}

\begin{center}
\uppercase{\bf Infinitely many absolute universes}
\vskip 20pt
{\bf Urban Larsson\footnote{larsson@iitb.ac.in}}\\
{\smallit Indian Institute of Technology Bombay, India}\\
{\bf Richard J.~Nowakowski\footnote{r.nowakowski@dal.ca}}\\
{\smallit Department of Mathematics and Statistics, Dalhousie University, Canada}\\
{\bf Carlos P. Santos\footnote{Under the scope of the projects
UIDB/00297/2020 and UIDP/00297/2020; cmf.santos@fc.ul.pt}}\\
%Under the scope of the projects UIDB/00297/2020 and UIDP/00297/2020;
%cmf.santos@fct.unl.pt
{\smallit Center for Mathematics and Applications (NovaMath), FCT NOVA, Portugal}\\
\end{center}

\begin{abstract}
Absolute combinatorial game theory was recently developed as a unifying tool for constructive/local game comparison (Larsson et al. 2018). The theory concerns {\em parental universes} of combinatorial games; standard closure properties are satisfied and each pair of non-empty sets of forms of the universe makes a form of the universe. Here we prove that there is an infinite number of absolute mis\`ere universes, by recursively expanding the dicot mis\`ere universe and the dead-ending universe. On the other hand, we prove that normal-play has exactly two absolute universes, namely the full space, and the universe of all-small games.
\end{abstract}

\section{Introduction}
%%%%%%%NEW INTRO%%%%%%%%%%%%%%%%%%
%\rjn{Classes of games that have the same properties, as in comparison tests, existence of inverses, \textit{etc}, reflect natural distinctions between types of rulesets. These can be determined by looking at `closure' properties. Under the Normal play convention, there are only two, the dicots and the full class of games. Under mis\`ere , both these classes are different. In this paper, we give a construction that yields  infinitely many different classes. A. Siegel \cite{Siegel} gives a different construction. These constructions arise naturally from actual games.}\\
Absolute combinatorial game theory was recently developed as a unifying tool for constructive/local game comparison \cite{LarSco, LarAbs}. A small but significant number of instances from the literature was used to motivate the theory. Here, we demonstrate that the absolute theory has infinitely many applications.

Recall the definition of partial order for short combinatorial games: $G\ge H$ if, for all games $X$, if Left wins $H+X$, then she wins $G+X$, as first and second player respectively. 
One of the most celebrated theorems in combinatorial game theory states that under the normal-play convention,  $G\ge H$ if and only if player Left wins $G-H$ playing second.

Under mis\`ere play, due to loss of group structure, comparison results are more intricate. But recent work find constructive methods in the sense that only the games $G$ and $H$, together with some proviso are required.  Since equivalence classes are small in full mis\`ere (for example the zero-class is trivial), one often turn towards restricted/modular mis\`ere, to find more order and algebraic structure. This naturally requires closure properties of standard game operators, such as taking options, addition and conjugation. Any set of games that satisfies the closure properties is a {\em universe} of games. 

In a recent break through the authors were able to prove a general constructive comparison result, given another natural closure property called {\em parentality}, that encompass not only normal and mis\`ere play, but also scoring play. A universe is parental if, for any pair of non-empty subsets of games $A$ and $B$ in the universe, the game $\{ A\mid B\}$ belongs to the universe. Universes that are parental and `outcome saturated' are called {\em absolute}, and the absolute theory applies. It turns out that outcome saturation is a consequence of parentality. Therefore the single relevant property to test whether a universe fits under the absolute umbrella is parentality.

Typically, a restriction of a convention concerns only the {\em game form}, and not the winning/ending convention. The most famous restriction is ``impartial'' which means that, at each stage of play, both players have the same options. However, the impartial universe is not parental. Another standard restriction is {\em dicot}, which means that, at each stage of play,  either both players or neither player can move. The universe of all dicots is parental (irrespective of winning/ending convention). Another restriction, with wide applications, is the universe of {\em dead-ending} games. If at some stage of play, player Left cannot move, then, after any sequence of Right moves, Left still cannot move. It is not possible that Right opens up any more Left options. And vice versa. All placement rulesets have this property. The dead-ending universe contains the dicotic universe, and it is parental. 

Siegel \cite{Si, Siege} study full mis\`ere, Dourbec et al. \cite{DRSS2015} study dicotic mis\`ere and Milley et al. study mis\`ere dead-ending \cite{LMNRS2021, MilleR2013}. Milley et al. \cite{MilleyRenault2018} give a brilliant survey on this development.  Analogous parental scoring play universes, that fit under the absolute umbrella, have been studied by Stewart (full scoring play) \cite{F1, F2}, Ettinger (dicot scoring play)\cite{Ettin2000}, and the authors (guaranteed scoring play)\cite{LarssNNS,LarSco}. 

Thus, we already know that absolute theory embraces several common situations. However, for a theory to be really useful, one would hope that it applies to an infinite set of instances. In this paper we prove that there is an infinite number of parental universes, with respect to distinguishability. We do this by recursively expanding the mis\`ere dicots and dead-ending games respectively. We prove that, in contrast, normal-play has only two absolute universes.

In Section~\ref{sec:abs} we recall some absolute terminology and basic combinatorial games' definitions. In Section~\ref{sec:parental} we define the notion of parentality. In Section~\ref{sec:closure} we discuss the dicot kernel of game forms and we define the (parental) closure of a given set of game forms. The rest of the paper is devoted to the classic conventions, Section~\ref{sec:classic}. We start with normal-play in Section~\ref{sec:normal}, and we finish off with the infinitely many mis\`ere absolute universes, Section~\ref{sec:misere}. In Section~\ref{sec:DE} we study universes between the dicot and  dead-ending universes, and in Section~\ref{sec:EOmega} we study universes that enclose the dead-ending universe.

\section{Some absolute terminology}\label{sec:abs}
The theory is built through a totally ordered, additive abelian, group of {\em adorns} $\A$.
A terminal position is of the form $\pura{\ell}{r}$ where $\ell, r\in \A$.  In general, if $G$ is a game with no Left options then we write $\GL=\emptyset^{\ell}$ for some $\ell\in \A$ and if Right has no options, we write $\GR=\emptyset^r$ for some $r\in \A$.

We refer to $\emptyset^a$ as an \textit{atom} and $a\in \A$ as the \textit{adorn}. A position in which Left (Right)
 does not have a move is called a \textit{Left- (Right-) atom}. A \textit{purely atomic} position
 is both Left- and Right-atomic.
Let $\bs{0}=\cg{\emptyset^0}{\emptyset^0}$, where $0$ is the identity of $\A$. In this paper, whenever $\A=\{0\}$, we will omit the adorn 0, and write simply $\emptyset=\atom{0}$.

A {\em free space} recursively builds all possible game forms, given a group of adorns.
 Let $\A$ be a group of adorns and let $\Omega_0=\{\pura{\ell}{r}\mid \ell, r \in \A\}$.
 For $n>0$, $\Omega_n$ is the set of all game forms with finite sets of options in $\Omega_{n-1}$,
 including game forms that are Left- and/or Right-atomic, and the set of game forms of \textit{birthday} $n$ is $\Omega_n\setminus \Omega_{n-1}$. Let $\Omega = \bigcup_{n\ge 0} \Omega_n$. Then $\Omega = (\Omega, \A)$ is a \textit{free space} of game forms.

By some abuse of terminology, we sometimes refer to the game form $G\in\Omega$ simply as a game, although there is not yet any incentive to play, or compare $G$ with other games in the same free space. An evaluation function $\nu: \A\times \A\rightarrow V$ takes care of the first part; here we are concerned exclusively with $A=0$ and $V=\{-1,+1\}$. The function $\nu$ distinguishes normal-play from mis\`ere: $\nu_L(0)=-1, \nu_R(0)=+1$ in normal-play, whereas  in mis\`ere we have $\nu_L(0)=+1, \nu_R(0)=-1$.

The notion of a {\em game space} $(\Omega,\A,V,\nu)$ includes the free space together with $E$ and $\nu$. Usually, we write $G\in \Omega$, and we assume an underlying given game space quadruple.

In a {\em disjunctive sum} of games, the current player picks exactly one game component and plays in it, by leaving any remaining components intact.
Here, and elsewhere, an expression of the type $\GL + H$ denotes the set of games  of the form $G^L+H$, $G^L\in \GL$, and this notion is only defined when $\GL$ is non-atomic.

Consider a free space $(\Omega,\A)$, and a pair of game forms $G, H\in (\Omega,\A)$. The disjunctive sum of $G$ and $H$ is given by:
\[G+H=
\begin{cases}
 \pura{\ell_1+\ell_2}{r_1+r_2} \textrm{ if } G=\pura{\ell_1}{r_1} \textrm{ and }
H=\pura{\ell_2}{r_2};\\
\cg{\emptyset^{\ell_1+\ell_2}}{\GR +H,G+\HR}, \textrm{ if }
G=\cg{\emptyset^{\ell_1}}{\GR},
H=\cg{\emptyset^{\ell_2}}{\HR},\\
 \textrm{ and at least one of $\GR $ and $\HR$
 non-atomic;}\\
\cg{\GL +H,G+\HL}{\emptyset^{r_1+r_2}}, \textrm{ if }
G=\cg{\GL}{\emptyset^{r_1}},
H=\cg{\HL}{\emptyset^{r_2}},\\
 \textrm{ and at least one of $\GL $ and $\HL$
 non-atomic;}\\
\cg{\GL +H,G+\HL}{\GR +H,G+\HR},
\textrm{ otherwise.}
\end{cases}\]

The \emph{conjugate} of a given game switches roles of the players.
The conjugate of $G\in \Omega$ is
\[ \conj{G} =
\begin{cases}
\pura{-r}{-\ell}, &\mbox{if $G=\pura{\ell}{r}$, $\ell,r\in \A$}\\
\cgs{\conj\GR}{\emptyset^{-\ell} }, &\mbox{if $G=\cgs{\emptyset^{\ell}}{\GR}$}\\
\cgs{\emptyset^{-r}}{\conj\GL }, &\mbox{if $G=\cgs{\GL}{\emptyset^{r}}$}\\
\cgs{\conj\GR }{\conj\GL }, &\mbox{otherwise},
\end{cases}
\]
where $\conj{\GL}$ denotes the set of game forms
$\conj{G^L} $, for $G^L\in \GL$, and similarly for $\GR$.

\subsection{Parental universes of game forms}\label{sec:parental}
A set of game forms $\U\subseteq(\Omega,\mathcal{A})$ is well behaved, and forms a {\em universe} of games, if it satisfies the following closure properties on the basic definitions.

\begin{enumerate}
  \item[a)] $\Omega_0\subseteq \U$;
  \item[b)] $\U$ is closed for disjunctive sum: if $G,H\in \U$, then $G+H\in \U$;
  \item[c)] $\U$ is closed for conjugation: if $G\in \U$, then $\conj{G}\in \U$;
  \item[d)] $\U$ is hereditary closed: if $G\in \U$ and $G'$ is an option of $G$, then $G'\in \U$;
\end{enumerate}

A set $S\subset (\Omega,\A)$ of game forms is {\em parental} if, for any pair of finite non-empty sets of game forms,
 $\cal G, \cal H\subset S$, then
 $\cg{\cal G}{\cal H}\in S$. If a parental set $S=\U$ is a universe, then $\U$ is a {\em parental universe}.\\

There are some classical parental universes based on the following definitions.

\begin{definition}[Dicot]\label{def:dicot}
A game form $G$ is a \emph{dicot} if $G^\mathcal{L}=\emptyset \Leftrightarrow G^\mathcal{R}=\emptyset$, and all options are dicots. The parental universe of dicotic forms is designated by $\D$.
\end{definition}

\begin{definition}[Dead-ending]\label{def:deadend}
A \emph{dead-ending} form has the property that if, at some stage of play, a player cannot
move, then they cannot move after any sequence of moves by the other player. A dead-ending form $G$ is called a \emph{Left-end} if $G^\mathcal{L}=\emptyset$ and a \emph{Right-end} if $G^\mathcal{R}=\emptyset$; thus, the zero
game $\{\emptyset\,|\,\emptyset\}$ is both a Left-end and a Right-end. A game is an \emph{end} if it is a Left-end or
a Right-end. The parental universe of dead-ending forms is designated by $\E$.
\end{definition}

\section{Closure of parental game forms}\label{sec:closure}
%This paper concerns parental universes.
\begin{problem}\label{prob:smallest}
For a given free space $\Omega$, what is the `smallest' parental universe?
\end{problem}
The notion of `smallest' requires a definition.
\begin{definition}[Kernel]\label{def:kernel}
Consider universes of games $ \U'\subseteq\U\subset (\Omega,\A)$, where $\U$ is parental. Then $\U$ is the kernel of $\Omega$, if it satisfies: if $\U'$ is a parental universe, then  $\U' = \U$.
\end{definition}

It turns out that the answer to Problem~\ref{prob:smallest} does not depend on the enclosing game space.

For any given game space, the smallest parental universe is the dicot universe.

\begin{theorem}[Dicot Forms]\label{thm:dicpar}
Consider a free space $(\Omega,\A)$, and  a universe $\U\subseteq(\Omega,\A)$. If $\U$ is parental then all dicot forms of $\Omega$ belong to $\U$.
\end{theorem}

\begin{proof}
Suppose that there is a dicot form $G\not\in\U$. Assume that $G$ is a simplest form in such conditions. We cannot have $G\in\Omega_0$ since, by definition of universe, $\Omega_0\subseteq \U$, contradicting the fact $G\not\in\U$. Therefore, we have $\GL\neq\emptyset$ and $\GR\neq\emptyset$. But, by definition of dicotic form, the elements of $\GL$ and $\GR$ are dicotic and, by the smallest rank assumption, those elements belong to $\U$. Therefore, by parentality, $G\in\U$; that is again a contradiction. All dicot forms of $\Omega$ belong to $\U$.
\end{proof}

\begin{corollary}[Dicot Kernel]\label{cor:dickern}
The kernel of a free space is its dicotic universe.
\end{corollary}
\begin{proof}
Combine Definitions~\ref{def:kernel} and~\ref{def:dicot} with Theorem~\ref{thm:dicpar}.
\end{proof}

Let us define three set operators on game forms.
\begin{definition}[Set Operators]
There are three set operators on a set of game forms $S\subseteq(\Omega,\A)$.
 \begin{enumerate}
   \item[] \hspace{-6 mm} The disjunctive sum operator: $\mathcal{D}(S)=\{G+H:G,H\in S\}$;
   \item[] \hspace{-6 mm} The conjugation operator: $\mathcal{C}(S)=\{\conj{G}:G\in S\}$;
   \item[] \hspace{-6 mm} The parental operator: $\mathcal{P}(S)=\{\cg{X}{Y}: X, Y\subseteq S, |X| > 0, |Y| > 0\}$.\\
 \end{enumerate}
\end{definition}

\begin{definition}[Closure of Game Forms]\label{def:Sclosure}
 Let $S\subseteq(\Omega,\A)$ be a hereditary closed set of game forms, and consider the following recursion:
\begin{enumerate}
\item[]\hspace{-6 mm}  Day $0$. $S_0=S\cup\Omega_0$;

\item[]\hspace{-6 mm}  Day $n>0$. $S_n=\mathcal{D}(\bigcup_{i=0}^{n-1} S_{i})\cup\mathcal{C}(\bigcup_{i=0}^{n-1} S_{i})\cup\mathcal{P}(\bigcup_{i=0}^{n-1} S_{i})$.
\end{enumerate}
The closure of $S$ is the set of game forms $\overline{S}=\bigcup_{i=0}^{\infty} S_{i}$.\\
\end{definition}

\begin{theorem}
Let $S\subseteq(\Omega,\A)$ be a hereditary closed set of game forms. The set $\overline{S}$ is a parental universe.
\end{theorem}
\begin{proof}
The hereditary closure follows by expanding Day 0, taking into account that $S$ is hereditary closed.
The closure for disjunctive sum, closure for conjugation, and parentality follows by Definition \ref{def:Sclosure}.
\end{proof}

\section{The classic conventions}\label{sec:classic}

\begin{problem}\label{prob:1}
How many parental universes are there in the classical conventions, normal- and mis\`ere-play respectively?
\end{problem}

We use the standard distinguishing operator for games.

Fix a convention, normal- or mis\`ere-play.
The perfect play outcome of a game $G$ is
\begin{itemize}
\item $o(G)=\Le$ if Left wins independently of who starts.
\item $o(G)=\Ri$ if Right wins independently of who starts.
\item $o(G)=\Ne$ if the next player wins independently of who starts.
\item $o(G)=\Pr$ if the previous player wins independently of who starts.
\end{itemize}
In the classic conventions the set of adorns is $\A=\{ 0\}$, and we may identify $\Omega = (\Omega, \A)$, and omit to adorn empty sets of games.

\begin{definition}[Distinguishability]
Let $\U\in\Omega$ be a universe. A pair of games $G,H\in\U$ are {\em distinguishable} modulo $\U$ if there is $X\in\U$ such that $o(G+X)\ne o(H+X)$. If $G$ and $H$ are indistinguishable modulo $\U$, we say they are \emph{equal} and write $G=H \pmod\U$.
\end{definition}

\begin{definition}[Proper Extension]\label{def:propext}
Let $\U,\U'\in\Omega$ be two universes such that $\U\subseteq\U'$. We say that $\U'$ is a \emph{proper extension} of $\U$ if there is $G\in\U'$, such that, for all $H\in\U$, $G$ and $H$ are distinguishable modulo $\U'$. If $\U'$ is not a proper extension of $\U$, we write $\U\simeq \U'$, meaning that the universes $\U$ and $\U'$ are the same with respect to distinguishability.
\end{definition}

\begin{observation} Note that universes may differ with respect to game forms, but be the same with respect to distinguishability.
\end{observation}

Returning to Problem \ref{prob:1}, the question that arises concerns how to obtain proper extensions from the Dicot Kernel. 

\subsection{Normal-play}\label{sec:normal}
Let us first analyze the game space $\Omega$ in normal-play.

\begin{theorem}\label{thm:normal}
Consider a normal-play parental universe $\U\subseteq\Omega$. If $\U$ contains a non-dicot form distinguished from zero modulo $\Omega$, then $\U\simeq\Omega$.
\end{theorem}

\begin{proof}
In Definition~\ref{def:propext}, take $\U'=\Omega$; all comparisons will concern the full game space of normal-play games.
We assume the existence of a non-dicot form in $\U$, and hence it has a follower $G$ that is Left-atomic or Right-atomic, but not both. Without loss of generality, suppose that $G$ is Right-atomic and not Left-atomic. Then Left wins playing second, so $G\su \bs 0$. But by assumption, $G\ne \bs 0$, so $G\succ \bs 0$.

We claim that $G=\bs k$ is a positive integer. Namely, $k=\min \{\ell : \bs \ell\succ G^L\}$, for all Left options of $G$, and where $\ell$ a positive integer.
Suppose first that Right starts the game $G-\bs k$. Since G is Right-atomic, the only option is to $G-\bs k +\bs 1$. By minimality of $k$, there is a left option $G^L$ such that $G^L\su \bs k-\bs 1$. Hence, Left wins when Right starts in $G^L-\bs k +\bs 1$.
Next, suppose that Left starts $G-\bs k$. By definition, for all Left options, $G^L-\bs k\prec 0$. Hence, Left loses.

We have the following possibilities:

\begin{enumerate}
  \item If $k=1$, then there is a form in $\U$ equal to $\bs 1$;
  \item If $k>1$, then, since by Simplicity Theorem $\{\bs 0\,|\,\bs k\}=\bs 1$ and due to the fact that $\U$ is parental, there is a form in $\U$ equal to $\bs 1$;
\end{enumerate}

In both cases, there is a form in $\U$ equal to $\bs 1$. But then, since $\U$ is closed for conjugation and for disjunctive sum, every integer has a form in $\U$.

Now, let $G=\{G\mathcal{^L}\,|\,G\mathcal{^R}\}$ be an arbitrary non-atomic canonical form under normal-play. Consider the form $G'$ obtained by replacing all the integers of $G\mathcal{^L}\cup G\mathcal{^R}$
by equal forms of $\U$, whose existence we have already guaranteed, and by replacing all the other forms of $G\mathcal{^L}\cup G\mathcal{^R}$ by equal forms of $\U$, whose existence is assumed by induction. On the one hand, by construction, $G'=_{\Np}G$ since the two forms have equivalent sets of options modulo $\Omega$. On the other hand, by parentally, $G'\in\U$. We conclude that each short game value has a form in $\U$, and hence $\U\simeq\Omega$.
\end{proof}

Thus, we have the following result.

\begin{corollary}
 The dicot (all-small) universe is the single proper absolute universe under normal-play.
\end{corollary}
\begin{proof}
 By Theorem~\ref{thm:normal}, we know that an absolute universe is no smaller than the universe of all dicots. By Theorem~\ref{thm:dicpar} we know that adjoining any non-dicot form to the dicots produces the full game space $\Omega$. Hence the universe of dicots $\U\subseteq \Omega$; note that if $X\in \U$, then $X+\bs 1\in\Le$, but $-\bs 1+\bs 1\in \Pr$, so indeed $\U\subset \Omega$ is a proper extension.
\end{proof}

\subsection{Mis\`ere-play}\label{sec:misere}
What about mis\`{e}re-play? Theorem~\ref{thm:mis} is the first step to the answer. We will give two different proofs, first a classical variation, by finding explicit distinguishing games, and then we use the modern approach, by absolute theory developed in \cite{LarAbs}.

We separate the outcome function in its two parts. The possible mis\`ere-play results are $\mathrm L$ (Left wins) and $\mathrm R$ (Right wins); by convention, they are totally ordered with $\mathrm L > \mathrm R$.
 The \emph{Left-} and \emph{Right-outcome}, in optimal play from both players, of a mis\`ere-play game $G$ is

 \[  o_L(G)= \left\{
\begin{array}{ll}
      {\mathrm L}, & \textrm{ if $\GL =\emptyset$;} \\
	{\max o_R(G^L)}, & \textrm{otherwise,} \\
\end{array}
\right. \]
\[o_R(G)= \left\{
\begin{array}{ll}
      {\mathrm R}, & \textrm{ if $\GR =\emptyset$;} \\
	{\min o_L(G^R)}, & \textrm{otherwise,}
\end{array}
\right. \]
respectively.

In \cite{LarAbs}, we proved the following result.

\begin{theorem}[Basic  Order, \cite{LarAbs}]\label{thm:basic}
Suppose $\U\subseteq(\Omega,\A )$ is an \abs universe of combinatorial games and let $G,H\in \U$. Then $G\succcurlyeq H$ if and only if the following two conditions hold\vspace{0.2 cm}

\noindent
Proviso:
\begin{enumerate}[]
 \item $o_L(G+X)\geqslant o_L(H+X)$ for all Left-atomic $X\in \U$;

 \item $o_R(G+X)\geqslant o_R(H+X)$ for all Right-atomic $X\in \U$;
\end{enumerate}
\noindent
Maintenance:
\begin{enumerate}[]
\item For all $G^R$, there is an $H^R$ such that $G^R\succcurlyeq H^R$, or there is a $G^{RL}$ such that $G^{RL}\succcurlyeq H$;
\item For all $H^L$, there is a $G^L$ such that $G^L\succcurlyeq H^L$, or there is an $H^{LR}$ such that $G\succcurlyeq H^{LR}$.
\end{enumerate}
\end{theorem}

A consequence of this is the following useful result. Its usefulness is due to that normal-play theory is well known, and so it can be easy to see when the normal-play inequality fails to hold. If it fails to hold, then the inequality fails in any (!) absolute universe. And indeed, in this paper we prove that there is an infinite number of such universes.

\begin{corollary}[Normal-play Order Preserving, \cite{LarAbs}]\label{thm:npop}
Let $G, H\in \U$, an \abs universe. If $G\succcurlyeq_\U H$ then $G\su_{\Np} H$.
\end{corollary}

The next results demonstrate that the mis\`ere-play convention has a richer structure than normal-play.

\subsubsection{Absolute universes strictly between $\D$ and $\E$}\label{sec:DE}

First, we present two proofs for the existence of an absolute universe strictly between $\D$ and $\E$. We will have use for a concept, closely related to the conjugate.

\begin{definition}[Adjoint]\label{def:adjoint}
Consider a mis\`ere-play universe $\U$ and a game $G\in \U$. Then $G^\circ$ is the adjoint of $G$, where
\[
G^\circ =
\begin{cases}
\cg{{\bs 0}}{{\bs 0}} & \textrm{if  $G = \bs 0$}; \\
\cg{\GR^\circ}{{\bs 0}} & \textrm{if $|\GR|>0$ and $\GL=\varnothing$}; \\
\cg{{\bs 0}}{\GL^\circ} & \textrm{if $|\GL|>0$ and $\GR=\varnothing$}; \\
\cg{\GR^\circ}{ \GL^\circ} & \textrm{otherwise},
\end{cases}
\]
and where $\circ$ applied to a set operates on its elements.
\end{definition}

It is well know that, under the mis\`ere-play convention, $G+G^\circ\in\mathcal{P}$.

\begin{theorem}\label{thm:mis}
Let $S=\D\cup\{\bs{1}\}$, where $\bs{1}=\{\bs 0\,|\,\varnothing \}$. Then, $\overline{S}$ is a proper extension of $\D$ and $\E$ is a proper extension of $\overline{S}$.
\end{theorem}

\begin{proof}[Proof 1.]
We observe that, by definition of closure of $S$, the forms $\bs 0$, $\bs{1}$, $\bs{1}+\bs{1}, \ldots$
are the only Right-atomic games in $\overline{S}$.

First, we prove that each Right-atomic game is distinguished from any dicot game. This follows by the Maintenance, specifically by Corollary~\ref{thm:npop}, because a strictly Right-atomic game is a positive number, and hence greater than any dicot in normal-play. Therefore, modulo $\overline{S}$, a Right-atomic game different than zero cannot be equal to a dicot. Hence $\overline{S}$ is a proper extension of $\D$.

Let us see that $\{\bs 0,\bs{1}\,|\,\varnothing\}$ differs from all Right-atomic games of $\overline{S}$  modulo $\E$.

\begin{itemize}
\item The game $\{\bs 0,\bs{1}\,|\,\varnothing\}\ne \bs 0$, because $o(\{\bs 0,\bs{1}\,|\,\varnothing\}) = \Ri$ and $o(\bs 0) = \Ne$.
\item The game $\{\bs 0,\bs{1}\,|\,\varnothing\}\ne \bs 1$, because $o(\{\bs 0,\bs{1}\,|\,\varnothing\}+\{\bs 0\,|\,*\}) = \Le$ and $o(\bs{1}+\{\bs 0\,|\,*\})= \Pr$.
\item The game $\{\bs 0,\bs{1}\,|\,\varnothing\}\ne \bs 1+\bs 1$, because $o(\{\bs 0,\bs{1}\,|\,\}+\{\bs 0\,|\,*\})= \Le$ and $o(\bs{1}+\bs{1}+\{\bs 0\,|\,*\}) = \Ne$.
\item The game $\{\bs 0,\bs{1}\,|\,\varnothing \}\neq \bs{1}+\bs{1}+\ldots+\bs 1$ (more than two summands) because $o(\{\bs 0,\bs{1}\,|\,\varnothing \}+\{\bs 0\,|\,*\}) = \Le$ and $o(\bs{1}+\bs{1}+\ldots+\bs 1+\{\bs 0\,|\,*\}) = \Ri$.
\end{itemize}

Now, let us consider a game $G$ not Right-atomic. Because $o(\{\bs 0,\bs{1}\,|\,\varnothing \})= \Ri$,  in case of equality, this forces $o(G)=\Ri$. Let $$X=\{ \{\bs 0\,|\,\text{all adjoints of followers of }G\}\,|\, *, \{\bs 0\,|\,*\}\}.$$
Playing first, Right wins $G+X$ by choosing the local mis\`{e}re-play winning line of $G$. This follows, since Left has only losing moves in $X$. Namely, whenever Left plays in $X$, Right responds to the appropriate adjoint-follower of $G$. In case of a final Left move to $\bs 0+X$, Right answers $*$; in case of a final Left move to $\bs{1}+X$, Right answers $\bs{1}+\{\bs 0\,|\,*\}$; in case of a final Left move to $\bs{1}+\ldots+\bs 1 + X$ (at least 2 summands), Right answers $\bs{1}+\ldots+\bs 1+*$.
On the other hand, playing first, Right loses $\{\bs 0,\bs{1}\,|\,\varnothing \}+X$.
Therefore, $\{\bs 0,\bs{1}\,|\,\varnothing \}\neq G$. This proves that $\E$ is a proper extension of $\overline{S}$.
\end{proof}

\begin{proof}[Proof 2.]
The first two paragraphs are the same as in the first proof.

Let us use Theorem~\ref{thm:basic} to show that $\{\bs 0,\bs{1}\,|\,\varnothing\}$ differs from all Right-atomic games of $\overline{S}$ modulo $\E$. We get extensive use of  Corollary~\ref{thm:npop}. In the third item, we cannot use it, but instead Maintenance gives a contradiction in $\E$, and moreover, the Proviso gives another contradiction.
\begin{itemize}
\item The game $\{\bs 0,\bs{1}\,|\,\varnothing\}\ne \bs 0$, because $\{\bs 0,\bs{1}\,|\,\varnothing\} =_\Np \bs{2}\succ_\Np \bs{0}$.
\item The game $\{\bs 0,\bs{1}\,|\,\varnothing\}\ne \bs 1$, because $\{\bs 0,\bs{1}\,|\,\varnothing\} =_\Np \bs{2}\succ_\Np \bs{1}$.
\item The game $\{\bs 0,\bs{1}\,|\,\varnothing\}\ne \bs 1+\bs 1$, because the Proviso is not satisfied with
 $\mathrm L=o_L(\{\bs 0,\bs{1}\,|\,\varnothing\}-\bs{1}) > o_L(\bs 1+\bs 1-\bs{1})=\mathrm R$. In fact, Maintenance is not satisfied, in spite the normal-play equality. Namely,  take $H=\{\bs 0,\bs{1}\,|\,\varnothing\}$ and $G=\bs 1+\bs 1$. Choose $H^L = \bs 0$. Then there is no $G^L$ such that $G^L\su H^L$, because $\bs 1\not\su \bs 0$ in mis\`ere play, and there is no $H^{LR}$ such that $G\su H^{LR}$.
\item The game $\{\bs 0,\bs{1}\,|\,\varnothing \}\neq \bs{1}+\bs{1}+\ldots+\bs 1$ (more than two summands) because $\{\bs 0,\bs{1}\,|\,\varnothing\} =_\Np 2\prec_\Np 3\pr \bs{1}+\bs{1}+\ldots+\bs 1$.
\end{itemize}

Moreover, by normal-play theory, and Corollary~\ref{thm:npop}, one may deduce that $\{\bs 0,\bs{1}\,|\varnothing\,\}\ne G$ whenever $\GR\neq\emptyset$ and $G\neq_\Np \bs{2}$. Otherwise, an argument similar to that used in the first proof is needed.
\end{proof}

Theorem~\ref{thm:mis} shows that the closure of $S=\D\cup\{\bs{1}\}$ is larger than the dicot universe but smaller than dead-ending universe. By extending this approach, we show that there is an infinite number of extensions strictly between the dicot and dead-ending universes.

\begin{definition}
Let $n$ be a nonnegative integer. Then, $n$ \emph{moves} for Left are recursively defined in the following way:
\begin{enumerate}
  \item[] $\bs{0}=\cg{\emptyset}{\emptyset}$;
  \item[] $\bs{n}=\{\bs {n}-\bs {1}\,|\,\emptyset\}$ if $n\geqslant 1$.
\end{enumerate}
The form $-\bs{n}$, $n$ \emph{moves} for Right, is recursively defined in a similar way.\\
\end{definition}

\begin{definition}
Let $n$ be a nonnegative integer. Then, $n$ \emph{controlled moves} for Left are recursively defined in the following way:
\begin{enumerate}
  \item[] $\widehat{0}=\bs{0}$;
  \item[] $\widehat{n}=\{\bs{0},\bs{1},\bs{2},\ldots,\bs{n}\bs{-1}\,|\,\emptyset\}$ if $n\geqslant 1$.
\end{enumerate}
The form $\widehat{-n}$, $n$ \emph{controlled moves} for Right, is recursively defined in a similar way.
\end{definition}

Notice that in normal-play, we don't need to distinguish moves from controlled moves. However, in mis\`ere-play, these forms are not equal. We will study the closure of the following universes.

\begin{definition}
 For $n\ge 1$, let $S_n=\D\cup \{\widehat{1},\ldots ,\widehat{n+1}\}$.
\end{definition}

Thus $\overline{S_0}$ is the universe studied in Theorem~\ref{thm:mis}. The universes $\overline{S_n}$ are parental, since they are extensions of the dicot universe.

\begin{definition}
Recursively, let $\circledast_n$ be the game form
\begin{enumerate}[]
  \item $\circledast_0=*$;
  \item $\circledast_{n}=\{0\,|\,\circledast_{n-1}\}$, if $n>0$.
\end{enumerate}
\end{definition}

\begin{observation}
Under the normal-play convention, these game forms are the canonical forms of the games $*,\uparrow,\Uparrow\!*, 3.\!\!\uparrow,4.\!\!\uparrow\!*,\ldots$
\end{observation}

\begin{lemma}\label{lem:Lwps}
Consider mis\`ere-play, with $n\ge 1$. Left wins $\bs{n}+\circledast_{k}$ playing second if and only if $n=k$.
\end{lemma}

\begin{proof}
If $k$ is zero, the result is trivial.

Suppose that $k>0$ and $n=k$. After a Right's first move from $\bs{k}+\circledast_{k}$ to $\bs{k}+\circledast_{k-1}$ then Left answers ${\bs k-\bs 1}+\circledast_{k-1}$ and wins by induction.

Suppose that $k> 0$ and $n\neq k$. After a Right's first move from $\bs{n}+\circledast_{k}$ to $\bs{n}+\circledast_{k-1}$, there are two possibilities. If Left replies from $\bs{n}+\circledast_{k-1}$ to $\bs{n}$, Right wins since he has no more moves.  If Left replies from $\bs{n}+\circledast_{k-1}$ to ${\bs n-\bs 1}+\circledast_{k-1}$, Right answers ${\bs n-\bs 1}+\circledast_{k-2}$ and, by induction, wins, or, if $k-1 = 0$, Right moves to ${\bs n - \bs 1}$ and also wins, since it leaves the last move to be made by Left. In this second case, observe that, since $k=1$,  $n\ge 1$, and $n\neq k$, we must have $n\geqslant 2$.
\end{proof}

The following result relies on equivalence classes in the mis\`{e}re-play convention. We give 2 proofs; in the first, we find an explicit distinguishing game, and in the second proof, we use absolute theory.

\begin{theorem}\label{thm:absinf}
There is an infinite number of absolute universes strictly between $\D$ and $\E$.
\end{theorem}

\begin{proof}[Proof 1.]
Let $n$ be a nonnegative integer. It suffices to prove that if $G\in\overline{S_n}$ then $G\neq\widehat{n+2}\mod{\overline{S_{n+1}}}$, that is, $\overline{S_{n+1}}$ is a proper extension of $\overline{S_{n}}$.

We observe that, by definition of closure of $S_n$, the forms $\widehat{k}$ with \linebreak $0\le k\leqslant n+1$, and all disjunctive sums of these forms, are the only Right-atomic games in $\overline{S_n}$. First, let us see that $\widehat{n+2}$ does not equal any Right-atomic game of $\overline{S_n}$.

Trivially, $\hat{n+2}\ne \bs 0$, since $o(\hat{n+2})=\Ri$ and $o(\bs 0)=\Ne$.

Let $X=\cgs{*}{\circledast_{0},\ldots, \circledast_{n+1}}$.  By Lemma~\ref{lem:Lwps}, Left wins playing second in $\widehat{n+2}+X$. Namely, Right must play to $\circledast{m}$, for some $0\le m\le n+1$, and then Left responds with the game $\hat{m}+\circledast{m}$.

Consider now an arbitrary Right-atomic disjunctive sum $G_1+\ldots+G_j\in \overline{S_n}$.

If this Right-atomic game has more than one summand, then Left loses $G_1+\ldots+G_j+X$ playing second, because Right goes to $G_1+\ldots+G_j+\circledast_0=G_1+\ldots+G_j+*$. Left loses because she has two or more moves in the Right-atomic summand.

If the disjunctive sum has only one summand, the best possible case for Left is $\widehat{n+1}+X$. (The ``hand tying principle'' holds in mis\`ere play for non-empty sets of options.) But, even so, Left loses playing second, because Right chooses $\widehat{n+1}+\circledast_{n+1}$ (Lemma~\ref{lem:Lwps}). Therefore, in all cases, Left loses playing second in $G_1+\ldots+G_j+X$.

We conclude that $\widehat{n+2}$ is not equal to a Right-atomic game of $\overline{S_n}$ modulo $\overline{S_{n+1}}$.

Suppose now that $o(G)=\Ri$, and let $$X=\{ \{0\,|\,\text{all adjoints of followers of }G)\}\,||\, \circledast_0,\circledast_1,\ldots,\circledast_{n+1}\}.$$
Playing first, Right wins $G+X$ by choosing the local mis\`{e}re-play winning line of $G$. By the adjoint construction, Left loses by playing in $X$, at any stage. In case of a Left move to $\bs 0+X$, Right answers $\circledast_0$; in case of a Left move to $G_1+\ldots+G_j+X$, with more than one Right-atomic summand, Right answers $G_1+\ldots+G_j+\circledast_0$; in case of a Left move to $G_1+X$, with one Right-atomic summand, Right answers $G_1+\circledast_{n+1}$. On the other hand, playing first, Right loses $\hat{n+2}+X$ (Lemma 15).
So, $\widehat{n+2}$ is different than $G$ modulo $\overline{S_{n+1}}$, and the proof is finished.
\end{proof}

\begin{proof}[Proof 2.]
It suffices to prove that if $G\in\overline{S_n}$ then $G\neq\widehat{n+2}\mod{\overline{S_{n+1}}}$.

We observe that, by definition of closure of $S_n$, the forms $\hat{k}$ with $0\le k\leqslant n+1$ and all disjunctive sums of these forms are the only Right-atomic games in $\overline{S_n}$. First, let us see that $\hat{n+2}$ does not equal any Right-atomic game $G$ of $\overline{S_n}$.

By Corollary~\ref{thm:npop}, $\hat{n+2}\ne  G$, unless possibly $G=\sum\hat{k_i}$, with $\sum k_i=n+2$, because $\hat{n+2}=_\Np n + 2$. Thus, it suffices to prove that Maintenance $\mod \overline{S_{n+1}}$ gives $G=\sum\hat{k_i}\ne \hat{n+2}=H$, whenever $\sum k_i = n + 2$. Since $\hat{n+2}\not\in \overline{S_n}$, the disjunctive sum has at least two summands. There is a Left option $H^L=\bs 0$, and for this option there is no $G^L$ such that $H^L\su G^L$. Namely any $G^L$ has at least one summand, $1\pr_{\Np} G^L$, and, so, that violates the normal-play inequality; by Corollary~\ref{thm:npop}, we know that $H^L=_{\Np} 0\prec_{\Np} 1\pr_{\Np} G^L$ cannot happen.

 Moreover, the option $H^L=\bs 0$, does not have a Right option, so we cannot have $H^{LR}\su G$.

We conclude that $\widehat{n+2}$ is not equal to a Right-atomic game of $\overline{S_n}$ modulo $\overline{S_{n+1}}$.

Also , $\hat{n+2}$ does not equal any $G$ with $G^\mathcal{R}\neq\emptyset$ modulo $\overline{S_{n+1}}$; the argument is the same as the proof of Theorem~\ref{thm:mis}.
\end{proof}

\subsubsection{Absolute universes strictly between $\E$ and $\Omega$}\label{sec:EOmega}

First, we present a proof for the existence of an absolute universe strictly between $\E$ and $\Omega$. From now on, for ease, we just present one proof.

\begin{theorem}\label{thm:mis2}
Let $S=\E\cup\{\{\varnothing |\,\bs{2}\}\}$. Then, $\overline{S}$ is a proper extension of $\E$ and $\Omega$ is a proper extension of $\overline{S}$.
\end{theorem}

\begin{proof}
We observe that, by definition of closure of $S$, the forms\linebreak $\{\varnothing |\,\bs{2}\}+\ldots+\{\varnothing |\,\bs{2}\}+L$, where $L$ is a Left-end form, are the only Left-atomic games of $\overline{S}$.

First, we prove that, modulo $\overline{S}$, $\{\varnothing |\,\bs{2}\}$ is distinguished from any Left-end form $G$. Let $X=\{-\bs{1}\,|\,\bs{0}\}$. This follows because, playing first,  Left wins $G+X$ by moving to $G-\bs{1}$; her victory is explained by the fact that $G$ is a Left-end. On the other hand, playing first, Left loses $\{\varnothing |\,\bs{2}\}+X$; she has to move to $\{\varnothing |\,\bs{2}\}-\bs{1}$ and Right replies $\bs{2}-\bs{1}$.

Second, let us consider a game $G\in\mathcal{N}$ not a Left-atomic form. Let\linebreak $X=\{\varnothing |\,\bs{2}\}+\ldots+\{\varnothing |\,\bs{2}\}$ with a sufficiently large number of copies of $\{\varnothing |\,\bs{2}\}$, depending on the rank of $G$. Left loses playing first in $G+X$; that happens because Right can give a large number of moves to Left, by opening the components after a first move by Left in $G$. On the other hand, Left wins playing first in $\{\varnothing |\,\bs{2}\}+X$; she has no moves. Hence, $\overline{S}$ is a proper extension of $\E$.

Third, consider the $\mathcal{N}$-position $\{\varnothing |\,\bs{3}\}\in\Omega$. We prove that, for any $G\in\overline{S}$, the form $\{\varnothing |\,\bs{3}\}$ is different than $G$ modulo $\Omega$. Of course, we only need to consider the case $G\in\mathcal{N}$. Let $$X=\{-\bs{2},\{\{-\bs{n}\,|\varnothing\}\,|\varnothing\}\ \mid \{\text{all adjoints of followers of }G\,|\,\bs{0}\}\},$$
with a sufficiently large $n$. Playing first, Left loses $\{\varnothing |\,\bs{3}\}+X$ since both $\{\varnothing |\,\bs{3}\}-\bs{2}$ and $\{\varnothing |\,\bs{3}\}+\{\{-\bs{n}\,|\varnothing\}\,|\varnothing \}$  are losing moves; Right replies to $\bs{3}-\bs{2}$ or to $\bs{3}+\{\{-\bs{n}\,|\varnothing \}\,|\varnothing\}$. On the other hand, Left wins $G+X$ playing first. Since $G$ is an $\mathcal{N}$-position, she can force a Right last move in the $G$ component. That must produce a position like $\{\varnothing |\,\bs{2}\}+\ldots+\{\varnothing |\,\bs{2}\}+L+X$ with Left to play. There are two possibilities:
\begin{enumerate}
  \item Right has two or more consecutive moves in $\{\varnothing |\bs{2}\}+\ldots+\{\varnothing |\,\bs{2}\}+L$; in that case, Left answers $\{\varnothing |\,\bs{2}\}+\ldots+\{\varnothing |\,\bs{2}\}+L+\{\{-\bs{n}\,|\varnothing \}\,|\varnothing\}$ and she is in time to reach $-\bs{n}$;
  \item $\{\varnothing |\,\bs{2}\}+\ldots+\{\varnothing |\,\bs{2}\}+L$ is $\{\varnothing |\,\bs{2}\}$ or $L$; in that case Left answers $\{\varnothing |\,\bs{2}\}-\bs{2}$ or $L-\bs{2}$ ($L$ is a Left-end).
\end{enumerate}
Hence, $\Omega$ is a proper extension of $\overline{S}$.
\end{proof}

\begin{observation} The reason why the extension is made with $\{\varnothing |\,\bs{2}\}$ and not with $\{\varnothing |\,\bs{1}\}$ is easy to understand. If Right plays on a disjunctive sum $\{\varnothing |\,\bs{2}\}+\ldots+\{\varnothing |\,\bs{2}\}$ with a large number of summands, he is able to ``give'' a large number of moves to Left, and that can be a big advantage under the mis\`ere-play convention. On the other hand, if Right plays on a disjunctive sum $\{\varnothing |\,\bs{1}\}+\ldots+\{\varnothing |\,\bs{1}\}$ with a large number of summands, he cannot ``give'' a large number of moves to Left, since whenever he opens a component, Left immediately finishes that component. In fact, $\overline{\E\cup\{\{\varnothing |\,\bs{1}\}\}}$ is not a proper extension of $\E$; using Theorem \ref{thm:basic}, it is possible to check that $\{\varnothing |\,\bs{1}\}$ is equal to $\bs{0}$ modulo $\overline{\E\cup\{\{\varnothing |\,\bs{1}\}\}}$.
\end{observation}

By extending the previous approach, we show that there is an infinite number of extensions strictly between $\E$ and $\Omega$. We need a preliminar lemma and a definition.

\begin{lemma}\label{thm:tip}
Let $G$ be the disjunctive sum $\{\varnothing |\,\bs{n_1}\}+\ldots+\{\varnothing |\,\bs{n_k}\}$ where $k\geqslant 2$ and where, for all $i$, $n_i\geqslant 2$. Let $n$ be a positive integer. If Right, playing first, wins $G-\bs{n}$, then Right, playing first, wins both
\begin{enumerate}
  \item $G-(\bs{n}+\bs{1})+*$, and
  \item $G+\{-(\bs{n}+\bs{1})+*\,|\,\bs 0\}$.
\end{enumerate}
\end{lemma}

\begin{proof}
Case 1. In the game $G-(\bs{n}+\bs{1})+*$, playing first, Right moves to $G-\bs{n}+*$. The only available answer for Left is to $G-\bs{n}$. Against that answer, by hypothesis, Right wins. Observe that, in total, Right makes $k+n+1$ moves before the $n_1+\ldots+n_k+1$ moves available for Left (the summand $1$ is the move in the star).

Case 2. In the game $G+\{-(\bs{n}+\bs{1})+*\,|\,\bs 0\}$, Right can move to $\bs{n_1}+\{\varnothing |\,\bs{n_2}\}+\ldots+\{\varnothing |\,\bs{n_k}\}+\{-(\bs{n}+\bs{1})+*\,|\,\bs 0\}$. A Left move to $(\bs{n_1}-\bs{1})+\{\varnothing |\,\bs{n_2}\}+\ldots+\{\varnothing |\,\bs{n_k}\}+\{-(\bs{n}+\bs{1})+*\,|\,\bs 0\}$ would lose the game since Right can answer with $(\bs{n_1}-\bs{1})+\{\varnothing |\,\bs{n_2}\}+\ldots+\{\varnothing |\,\bs{n_k}\}$, opening all the remaining components after a move that Left still has in the first one. Hence, let us analyse a Left move to $\bs{n_1}+\{\varnothing |\,\bs{n_2}\}+\ldots+\{\varnothing |\,\bs{n_k}\}-(\bs{n}+\bs{1})+*$. Against that move, Right adopts the strategy of playing consecutively in the component $-(\bs{n}+\bs{1})$. 

Two things may happen: (1) Right runs out of moves in that component before Left has to play in the star, he opens the second summand of $G$ and wins; (2) Left spends all her moves in the components $\bs{n_1}$ and $*$ before Right runs out of moves in $-(\bs{n}+\bs{1})$, Right makes all the moves as in the analysis of $G-(\bs{n}+\bs{1})+*$. In total, Right makes $k+n+1$ moves before the $n_1+\ldots+n_k+2$ moves for Left (the summand $2$ respects to the two Left moves in $\{-(\bs{n}+\bs{1})+*\,|\,\bs 0\}$ and in $*$). This is even one more move for Left than in the case of the analysis of the game $G-(\bs{n}+\bs{1})+*$. Therefore, Right also wins in this second case.
\end{proof}

\begin{definition}
Let $n\geqslant 2$ be an integer. Then, the \emph{Left-hook} of order $n$ is the game form $\zs{n}=\{\varnothing |\,\bs{n}\}$. The \emph{Right-hook} of order $n$ is the game form $\zs{-n}=\{-\bs{n}\,|\varnothing \}$.
\end{definition}

\begin{definition}\label{def:hook}
 For $n\ge 2$, let $Z_n=\E\cup \{\zs{2},\ldots ,\zs{n}\}$.
\end{definition}

\begin{theorem}\label{thm:absinf2}
There is an infinite number of absolute universes strictly between $\E$ and $\Omega$.
\end{theorem}

\begin{proof}
It suffices to prove that, for all $n\geqslant 2$, if $G\in\overline{Z_n}$ then $G\neq\zs{n+1}\mod{\overline{Z_{n+1}}}$. Indeed, this proves that $\overline{Z_{n+1}}$ is a proper extension of $\overline{Z_{n}}$.

We observe that, by the definition of closure of $Z_n$ (Definition \ref{def:Sclosure}), the forms $\zs{n_1}+\ldots+\zs{n_k}+L$, with  $2\le n_i\leqslant n$ and a Left-end $L$, are the only Left-atomic games in $\overline{Z_n}$. First, let us see that $\zs{n+1}$ does not equal any Left-atomic game of $\overline{Z_n}$.

If $L\not\cong \{\varnothing |\varnothing \}$ then $\zs{n_1}+\ldots+\zs{n_k}+L$ is a negative integer under the normal-play convention and $\zs{n+1}$ is equal to zero under the normal-play convention. Due to that, by Corollary \ref{thm:npop}, the games $\zs{n_1}+\ldots+\zs{n_k}+L$ and $\zs{n+1}$ are not equal modulo $\overline{Z_{n+1}}$. Hence, suppose that $L\cong \{\varnothing |\varnothing \}$, that is, consider only disjunctive sums of the type $\zs{n_1}+\ldots+\zs{n_k}$.

If there is only one summand, let $X=-\bs{n}$. We have that, playing first, Right wins $\zs{n+1}+X$ by moving to $\bs{n}+\bs{1}-\bs{n}$, and loses $\zs{n_1}+X$ since $n_1\leqslant n$. Hence, suppose that there are two or more summands.

Let $X=-\bs{n}$. If, playing first, Right loses $\zs{n_1}+\ldots+\zs{n_k}+X$, then, since Right wins $\zs{n+1}+X$, $\zs{n_1}+\ldots+\zs{n_k}$ and $\zs{n+1}$ are not equal modulo $\overline{Z_{n+1}}$. Therefore, suppose that, playing first, Right wins $\zs{n_1}+\ldots+\zs{n_k}-\bs{n}$. By Lemma~\ref{thm:tip}, Right also wins $\zs{n_1}+\ldots+\zs{n_k}+\{-(\bs{n}+\bs{1})+*\,|\,\bs 0\}$. But, playing first, Right loses $\zs{n+1}+\{-(\bs{n}+\bs{1})+*\,|\,\bs 0\}$; if he moves to $\zs{n+1}$, Left runs out of moves and wins; if he moves to $(\bs{n}+\bs{1})+\{-(\bs{n}+\bs{1})+*\,|\,\bs 0\}$, Left wins with the move $(\bs{n}+\bs{1})-(\bs{n}+\bs{1})+*$ since she can force an extra Right's move in the star. Therefore, in this last case, $X=\{-(\bs{n}+\bs{1})+*\,|\,\bs 0\}$ distinguishes $\zs{n_1}+\ldots+\zs{n_k}$ from $\zs{n+1}$.

Suppose now that $G\in\mathcal{N}$ has a Left option. Let $X=\zs{2}+\ldots+\zs{2}$ be a disjunctive sum with a sufficiently large number of copies of $\zs{2}$, depending on the rank of $G$. Left loses playing first in $G+X$; that happens because Right is in time to give a very large number of moves to Left, by opening the components after a first move by Left in $G$. On the other hand, Left wins playing first in $\zs{n+1}+X$; she has no moves. $\overline{Z_{n+1}}$ is a proper extension of $\overline{Z_{n}}$.
\end{proof}

%\clearpage

\end{document}